\overfullrule=0pt
\centerline {\bf Energy functionals of Kirchhoff-type problems having multiple global minima}\par
\bigskip
\bigskip
\centerline {BIAGIO RICCERI}\par
\bigskip
\bigskip
{\bf Abstract.} In this paper, using the theory developed in [8], we obtain some results of a totally new type about a class of
non-local problems. Here is a sample:\par Let $\Omega\subset {\bf R}^n$ be a smooth bounded domain, with $n\geq 4$, let $a, b, \nu\in {\bf R}$,
with $a\geq 0$ and $b>0$, and let 
 $p\in \left ] 0,{{n+2}\over {n-2}}\right [$.\par
Then, for each $\lambda>0$ large enough and for each convex set $C\subseteq L^2(\Omega)$ whose closure in $L^2(\Omega)$
contains $H^1_0(\Omega)$, there exists $v^*\in C$ such that the problem
$$\cases {-\left ( a+b\int_{\Omega}|\nabla u(x)|^2dx\right )\Delta u =\nu|u|^{p-1}u+\lambda(u-v^*(x)) & in $\Omega$\cr
& \cr u=0 & on $\partial\Omega$\cr}$$
has at least three weak solutions, two of which are global minima in $H^1_0(\Omega)$ of the corresponding energy functional.\par
\bigskip
\bigskip
\bigskip
\bigskip
{\bf Introduction}\par
\bigskip
\bigskip
Here and in what follows, $\Omega\subset {\bf R}^n$ ($n\geq 3$) is a bounded domain with smooth boundary and $a, b\in {\bf R}$,
with $a\geq 0$ and $b>0$.\par
\par
Consider the non-local problem
$$\cases {-\left ( a+b\int_{\Omega}|\nabla u(x)|^2dx\right )\Delta u =h(x,u) & in $\Omega$\cr
& \cr u=0 & on $\partial\Omega$\ ,\cr}$$
$h:\Omega\times {\bf R}\to {\bf R}$ being a Carath\'eodory function.\par
\smallskip
In the past years, many papers have been produced on the existence and multiplicity of weak solutions for this problem.
Usual reference papers are [1], [2], [4], [5], [6], [10].
\smallskip
In the present note, we are interested in the multiplicity of solutions of the above problem under the following new aspect which has
not been considered before: among the solutions of the problem, at least two are global minima of the energy functional.\par
\smallskip
But, more generally, the global structure itself of the conclusions we reach is novel at all, so that no proper comparison of our results with the previous ones 
in the field can be made.
\par
\smallskip
For instance, the following proposition summarizes very well the novelties of our main result (Theorem 1) of which it is the simplest particular case:\par
\medskip
PROPOSITION 1. - {\it Let $n\geq 4$, let $\nu\in {\bf R}$ and let $p\in \left ] 0,{{n+2}\over {n-2}}\right [$.\par
Then, for each $\lambda>0$ large enough and for each convex set $C\subseteq L^2(\Omega)$ whose closure in $L^2(\Omega)$
contains $H^1_0(\Omega)$, there exists $v^*\in C$ such that the problem
$$\cases {-\left ( a+b\int_{\Omega}|\nabla u(x)|^2dx\right )\Delta u =\nu|u|^{p-1}u+\lambda(u-v^*(x)) & in $\Omega$\cr
& \cr u=0 & on $\partial\Omega$\cr}$$
has at least three solutions, two of which are global minima in $H^1_0(\Omega)$ of the functional
$$u\to {{a}\over {2}}\int_{\Omega}|\nabla u(x)|^2dx+{{b}\over {4}}\left ( \int_{\Omega}|\nabla u(x)|^2dx\right ) ^2-{{\nu}\over
{p+1}}\int_{\Omega}|u(x)|^{p+1}dx
-{{\lambda}\over {2}}\int_{\Omega}|u(x)-v^*(x)|^2dx\ .$$}
\medskip
What allows us to obtain results of this totally new type is the use of theory we have recently developed in [8].
\bigskip
{\bf Results}\par
\bigskip
On the Sobolev space $H^1_0(\Omega)$, we consider the scalar product
$$\langle u,v\rangle=\int_{\Omega}\nabla u(x)\nabla v(x)dx$$
and the induced norm
$$\|u\|=\left ( \int_{\Omega}|\nabla u(x)|^2dx\right )^{1\over 2}\ .$$
We denote by ${\cal A}$ the class of all Carath\'eodory functions $f:\Omega\times {\bf R}\to {\bf R}$ such that
$$\sup_{(x,\xi)\in\Omega\times {\bf R}}{{|f(x,\xi)|}\over {1+|\xi|^p}}<+\infty \eqno{(1)}$$
for some $p\in \left ]0, {{n+2}\over {n-2}}\right [$.\par
\smallskip
Moreover, we denote by $\tilde{\cal A}$ the class of all Carath\'eodory functions $g:\Omega\times {\bf R}\to {\bf R}$ such that
$$\sup_{(x,\xi)\in\Omega\times {\bf R}}{{|g(x,\xi)|}\over {1+|\xi|^q}}<+\infty \eqno{(2)}$$
for some $q\in \left ] 0, {{2}\over {n-2}}\right [$.
\smallskip
Furthermore, we denote by $\hat {\cal A}$ the class of all functions $h:\Omega\times {\bf R}\to {\bf R}$ of the type
$$h(x,\xi)=f(x,\xi)+\alpha(x)g(x,\xi)$$
with $f\in {\cal A}$, $g\in\tilde{\cal A}$ and $\alpha\in L^2(\Omega)$.
For each $h\in\hat{\cal A}$, we define the functional $I_h:H^1_0(\Omega)\to {\bf R}$, by putting
$$I_h(u)=\int_{\Omega}H(x,u(x))dx$$
for all $u\in H^1_0(\Omega)$, where
$$H(x,\xi)=\int_0^{\xi}h(x,t)dt$$
for all $(x,\xi)\in \Omega\times {\bf R}$.\par
\smallskip
By classical results (involving the Sobolev embedding theorem), the functional $I_h$
turns out to be sequentially weakly continuous, of class $C^1$, with compact derivative given by
$$I_h'(u)(w)=\int_{\Omega}h(x,u(x))w(x)dx$$
for all $u,w\in H^1_0(\Omega)$.\par
\smallskip
Now, let us recall that, given $h\in\hat{\cal A}$, a weak solution of the problem
$$\cases {-\left ( a+b\int_{\Omega}|\nabla u(x)|^2dx\right )\Delta u =h(x,u) & in $\Omega$\cr
& \cr u=0 & on $\partial\Omega$\cr}$$
is any $u\in H^1_0(\Omega)$ such that
$$\left ( a+b\int_{\Omega}|\nabla u(x)|^2dx\right )\int_{\Omega}\nabla u(x)\nabla w(x)dx=
\int_{\Omega}h(x,u(x))w(x)dx$$
for all $w\in H^1_0(\Omega)$. Let $\Phi:H^1_0(\Omega)\to {\bf R}$ be the functional defined by
$$\Phi(u)={{a}\over {2}}\|u\|^2+{{b}\over {4}}\|u\|^4$$
for all $u\in H^1_0(\Omega)$.\par
\smallskip
Hence, the weak solutions of the problem are precisely the critical points in
$H^1_0(\Omega)$ of the functional $\Phi-I_h$ which is said to be the energy functional of the problem.\par
\smallskip
Here is our main result:\par
\medskip
THEOREM 1. - {\it Let $n\geq 4$, let $f\in {\cal A}$ and let $g\in\tilde{\cal A}$ be such
that the set
$$\left \{x\in \Omega : \sup_{\xi\in {\bf R}}|g(x,\xi)|>0\right\}$$
has a positive measure.\par
Then,  there exist $\lambda^*\geq 0$ such that, for each $\lambda>\lambda^*$ and for each 
convex set $C\subseteq L^2(\Omega)$ whose closure in $L^2(\Omega)$
contains the set $\{G(\cdot,u(\cdot)) : u\in H^1_0(\Omega)\}$,
there exists $v^*\in C$ such
that the problem
$$\cases {-\left ( a+b\int_{\Omega}|\nabla u(x)|^2dx\right )\Delta u =f(x,u)+\lambda(G(x,u)-v^*(x))g(x,u) & in $\Omega$\cr
& \cr u=0 & on $\partial\Omega$\cr}$$
has at least three weak solutions, two of which are global minima in $H^1_0(\Omega)$ of the functional
$$u\to {{a}\over {2}}\int_{\Omega}|\nabla u(x)|^2dx+{{b}\over {4}}\left ( \int_{\Omega}|\nabla u(x)|^2dx\right ) ^2-\int_{\Omega}F(x,u(x))dx
-{{\lambda}\over {2}}\int_{\Omega}|G(x,u(x))-v^*(x)|^2dx\ .$$
If, in addition,  the functional
$$u\to {{a}\over {2}}\int_{\Omega}|\nabla u(x)|^2dx+{{b}\over {4}}\left ( \int_{\Omega}|\nabla u(x)|^2dx\right ) ^2-\int_{\Omega}F(x,u(x))dx$$
has at least two global minima in $H^1_0(\Omega)$ and the function $G(x,\cdot)$ is strictly monotone for all $x\in \Omega$,
then $\lambda^*=0$.}\par
\medskip
The main tool we use to prove Theorem 1 is Theorem C below which is, in turn, a direct consequence of
two other results recently established in [8].\par
\medskip
To state Theorem C in a compact form, we now introduce some notations.\par
\smallskip
Here and in what follows, $X$ is a non-empty set, $V, Y$ are two topological
spaces, $y_0$ is a point in $Y$.\par
\smallskip
We denote by ${\cal G}$
the family of all lower semicontinuous functions $\varphi:Y\to [0,+\infty[$,
with $\varphi^{-1}(0)=\{y_0\}$, such that, for each neighbourhood $U$ of $y_0$,
one has
$$\inf_{Y\setminus U}\varphi>0\ .$$
Moreover, we denote by ${\cal H}$ the family of all functions
$\Psi:X\times V\to Y$ such that, for each $x\in X$, 
$\Psi(x,\cdot)$ is continuous, injective, open, 
takes the value $y_0$ at a point $v_x$ and the function
$x\to v_x$ is not constant. 
Furthermore, we denote by ${\cal M}$ the family of all
functions $J:X\to {\bf R}$ whose set of all global minima
(noted by $M_{J}$) is non-empty.\par
\smallskip
Finally, for each $\varphi\in {\cal G}$, $\Psi\in {\cal H}$ 
 and $J\in {\cal M}$, we put
$$\theta(\varphi,\Psi,J)=\inf\left \{
{{J(x)-J(u)}\over {\varphi(\Psi(x,v_u))}} : (u,x)\in M_{J}\times
X\hskip 3pt \hbox {\rm with}\hskip 3pt v_x\neq v_u\right
\}\ .$$
When $X$ is a topological space, a function $\psi:X\to {\bf R}$ is said to be 
inf-compact if $\psi^{-1}(]-\infty,r])$ is compact for all $r\in {\bf R}$.
\medskip
THEOREM A ([8], Theorem 3.1). - {\it Let $\varphi\in {\cal G}$, $\Psi\in {\cal H}$  
 and $J\in {\cal M}$. \par
Then, for each $\lambda>\theta(\varphi,\Psi,J)$, one has
$$\sup_{v\in V}\inf_{x\in X}
(J(x)-\lambda\varphi(\Psi(x,v)))<
\inf_{x\in X}\sup_{z\in X}
(J(x)-\lambda\varphi(\Psi(x,v_z)))\ .$$}\par
\medskip
THEOREM B ([8], Theorem 3.2). - {\it Let $X$ be a topological space, $E$ a real Hausdorff
topological vector space, $C\subseteq
E$ a convex set, $f:X\times C\to
{\bf R}$ a function which is lower semicontinuous and inf-compact in $X$,
and upper semicontinuous and concave in $C$. Assume also that
$$\sup_{v\in C}\inf_{x\in X}f(x,v)<
\inf_{x\in X}\sup_{v\in C}f(x,v)\ .$$
Then, there exists $v^*\in C$ such that
the function $f(\cdot,v^*)$ has at least two global minima in $X$.}\par
\medskip
THEOREM C. - {\it Let $\varphi\in {\cal G}$, $\Psi\in {\cal H}$  
 and $J\in {\cal M}$. 
Moreover, assume that $X$ is a topological space, that
$V$ is a real Hausdorff topological vector space and that
$\varphi(\Psi(x,\cdot))$ is convex and continuous for each $x\in X$. Finally, let
$\lambda>\theta(\varphi,\Psi,J)$  and let $C\subseteq V$ be a convex set, with
$\{v_x : x\in X\}\subseteq \overline {C}$,
such that
 the function $x\to J(x)-\lambda\varphi(\Psi(x,v))$ is
lower semicontinuous and inf-compact in $X$ for all $v\in C$.\par
Under such hypotheses, 
there exists $v^*\in C$
 such
that the function $x\to J(x)-\lambda\varphi(\Psi(x,v^*))$ 
has at least two global minima in $X$.}\par
\smallskip
PROOF. Set
$$D=\{v_x : x\in X\}$$
and, for each $(x,v)\in X\times V$, put
$$f(x,v)=J(x)-\lambda\varphi(\Psi(x,v))\ .$$
Theorem A ensures that
$$\sup_{v\in V}\inf_{x\in X}f(x,v)<
\inf_{x\in X}\sup_{v\in D}f(x,v)\ .\eqno{(3)}$$
But, since $f(x,\cdot)$ is continuous and $D\subseteq \overline {C}$,
 we have
$$\sup_{v\in D}f(x,v)=\sup_{v\in \overline {D}}f(x,v)\leq \sup_{v\in \overline {C}}f(x,v)
= \sup_{v\in C}f(x,v)$$
for all $x\in X$, and hence, taking $(3)$ into account, it follows that 
$$\sup_{v\in C}\inf_{x\in X}f(x,v)<
\inf_{x\in X}\sup_{v\in D}f(x,v)\leq \inf_{x\in X}\sup_{v\in C}f(x,v)\ .$$
At this point, the conclusion follows applying Theorem B to the restriction of
the function $f$ to $X\times C$.\hfill $\bigtriangleup$\par
\medskip
{\it Proof of Theorem 1.} For each $\lambda\geq 0$,
$v\in L^2(\Omega)$, consider the function $h_{\lambda,v}:\Omega\times {\bf R}\to {\bf R}$
defined by 
$$h_{\lambda,v}(x,\xi)=f(x,\xi)+\lambda(G(x,\xi)-v(x))g(x,\xi)$$ 
for all $(x,\xi)\in\Omega\times {\bf R}$.
Clearly, the function $h_{\lambda,v}$ lies in
$\hat{\cal A}$ and
$$H_{\lambda,v}(x,\xi)=F(x,\xi)+{{\lambda}\over {2}}\left ( |G(x,\xi)-v(x)|^2-|v(x)|^2\right )\ .$$
Moreover, if $p\in \left] 0,{{n+2}\over {n-2}}\right [$ and $q\in
\left ] 0,{{2}\over {n-2}}\right [$ are such that $(1)$ and $(2)$ hold, 
 for some constant $c_{\lambda,v}$, we have
$$\int_{\Omega}|H_{\lambda,v}(x,u(x))|dx\leq c_{\lambda,v}\left ( \int_{\Omega}|u(x)|^{p+1}dx+
\int_{\Omega}|u(x)|^{2(q+1)}dx+1\right )$$
for all $u\in H^1_0(\Omega)$. Therefore, by the Sobolev embedding theorem, for a constant $\tilde c_{\lambda,v}$,
we have
$$\Phi(u)-I_{h_{\lambda,v}}(u)\geq {{b}\over {4}}\|u\|^4-\tilde c_{\lambda,v}(\|u\|^{p+1}+\|u\|^{2(q+1)}+1) \eqno{(4)}$$
for all $u\in H^1_0(\Omega)$. On the other hand, since $n\geq 4$, one has
$$\max\{p+1, 2(q+1)\}<{{2n}\over {n-2}}\leq 4\ .$$
Consequently, from $(4)$, we infer that
$$\lim_{\|u\|\to +\infty}(\Phi(u)-I_{h_{\lambda,v}}(u))=+\infty\ .\eqno{(5)}$$
Since the functional $\Phi-I_{h_{\lambda,v}}$ is sequentially weakly lower semicontinuous, by the Eberlein-Smulyan theorem and
by $(5)$, it follows that it is inf-weakly compact.\par
Now, we are going to apply Theorem C taking $X=H^1_0(\Omega)$ with
the weak topology, $V=Y=L^2(\Omega)$ with the strong topology and $y_0=0$. Also, we take
$$\varphi(w)={{1}\over {2}}\int_{\Omega}|w(x)|^2dx$$
for all $w\in L^2(\Omega)$. Clearly, $\varphi\in {\cal G}$. Furthermore, we take
$$\Psi(u,v)(x)=G(x,u(x))-v(x)$$
for all $u\in H^1_0(\Omega)$, $v\in L^2(\Omega)$, $x\in \Omega$. Clearly, $\Psi(u,v)\in L^2(\Omega)$,
 $\Psi(u,\cdot)$ is a homeomorphism, and we have
$$v_u(x)=G(x,u(x))\ .$$
We now show that the map $u\to v_u$ is not constant in $H^1_0(\Omega)$. Set
$$A=\left \{x\in \Omega : \sup_{\xi\in {\bf R}}|g(x,\xi)|>0\right\}\ .$$
By assumption, meas$(A)>0$. Then, by the classical Scorza-Dragoni theorem ([3], Theorem 2.5.19),
there exists a compact set $K\subset A$, of positive measure, such that the restriction
of $G$ to $K\times {\bf R}$ is continuous. Fix a point $\tilde x\in K$ such that
the intersection of $K$ and any ball centered at $\tilde x$ has a positive measure. Of course, the function
$G(\tilde x,\cdot)$ is not constant.
Fix $\xi_1, \xi_2\in {\bf R}$
such that
$$G(\tilde x,\xi_1)<G(\tilde x,\xi_2)\ .$$
By continuity, there is a closed ball $B(\tilde x,r)$ such that
$$G(x,\xi_1)<G(x,\xi_2)$$
for all $x\in K\cap B(\tilde x,r)$. Finally, consider two functions $u_1, u_2\in H^1_0(\Omega)$
such that
$$u_1(x)=\xi_1$$
and
$$u_2(x)=\xi_2$$
for all $x\in K\cap B(\tilde x,r)$. So
$$G(x,u_1(x))<G(x,u_2(x))$$
for all $x\in K\cap B(\tilde x,r)$. Hence, as meas($K\cap B(\tilde x,r))>0$, we infer that
$v_{u_1}\neq v_{u_2}$, as claimed. As a consequence, $\Psi\in {\cal H}$.  Of course,
$\varphi(\Psi(u,\cdot))$  is continuous and convex for all $u\in X$. 
Finally, take
$$J=\Phi-I_f\ .$$
Clearly, $J\in {\cal M}$. So,  for what seen above, all the assumptions of Theorem C are satisfied.
Consequently, if we take 
$$\lambda^*=\theta(\varphi, \Psi,J)\eqno{(6)}$$ and
fix $\lambda>\lambda^*$ and a convex set  $C\subseteq L^2(\Omega)$ whose closure in $L^2(\Omega)$
contains the set $\{G(\cdot,u(\cdot)) : u\in H^1_0(\Omega)\}$, there exists $v^*\in C$ such that
the functional $\Phi-I_{h_{\lambda,v^*}}$ has at least two global minima in $H^1_0(\Omega)$ which are, therefore,
weak solutions of the problem we are dealing with. To guarantee the existence of a third solution, 
  denote by $k$ the inverse of the restriction of the function $at+bt^3$ to $[0,+\infty[$. 
Let $T:X\to X$ be the operator defined by
$$T(w)=\cases {{{k(\|w\|)}\over {\|w\|}}w & if $w\neq 0$\cr & \cr
0 & if $w=0$\ ,\cr}$$
 Since $k$ is continuous
and $k(0)=0$, the operator $T$ is continuous in $X$. For each $u\in X\setminus \{0\}$,
we have
$$T(\Phi'(u))=T((a+b\|u\|^2)u)={{k((a+b\|u\|^2)\|u\|)}\over {(a+b\|u\|^2)\|u\|}}(a+b\|u\|^2)u={{\|u\|}\over {(a+b\|u\|^2)\|u\|}}(a+b\|u\|^2)u=u\ .$$
In other words, $T$ is a continuous inverse of $\Phi'$. Then, since $I_{h_{\lambda,v^*}}'$ is compact, the functional
$\Phi-I_{h_{\lambda,v^*}}$ satisfies the Palais-Smale condition ([9], Example 38.25) and hence the existence of a third critical
point of the same functional is assured by Corollary 1 of [7].\par
Finally, assume that the functional $\Phi-I_f$ has at least two global minima, say $\hat u_1, \hat u_2$. So,
the set 
$$\Gamma:=\{x\in \Omega : \hat u_1(x)\neq \hat u_2(x)\}$$ 
has a positive measure. By assumption, we have
$$G(x,\hat u_1(x))\neq G(x,\hat u_2(x))$$
for all $x\in \Gamma$, and hence $v_{\hat u_1}\neq v_{\hat u_2}$. Then, by definition, we have
$$0\leq\theta(\varphi,\Psi,J)\leq {{J(\hat u_1)-J(\hat u_2)}\over {\varphi(\Psi(\hat u_1,v_{\hat u_2}))}}=0$$
and so $\lambda^*=0$ in view of $(6)$.\hfill $\bigtriangleup$\par
\medskip
REMARK 1. - In Theorem 1,  the assumption made on $g$ (besides $g\in \tilde{\cal A}$) is essential. Indeed, if $g(x,\xi)=0$
for almost every $x\in \Omega$ and for every $\xi\in {\bf R}$, then, taking
$f=0$ (which is an allowed choice), the problem would have the zero solution only.\par
\medskip
REMARK 2. - The assumption $n\geq 4$ is likewise essential. Indeed, Proposition 1 (which comes from Theorem 1 taking
$f(x,\xi)=\nu|\xi|^{p-1}\xi$ and $g(x,\xi)=1$) does not hold if $n=3$ and $\nu>0$. To see this, 
 take $p=4$ (which is compatible
with the condition $p<{{n+2}\over {n-2}}$ when $n=3$) and observe
that the corresponding energy functional is unbounded below.\par
\medskip
Besides Proposition 1, among the consequences of Theorem 1, we highlight the following\par
\medskip
THEOREM 2. - {\it Let $n\geq 4$, let $f\in {\cal A}$ and let $g\in\tilde{\cal A}$ be such
the set
$$\left\{x\in \Omega : \sup_{\xi\in {\bf R}}F(x,\xi)>0\right\}$$
has a positive measure.
 Moreover, assume that, for each $x\in \Omega$, $f(x,\cdot)$ is
odd, $g(x,\cdot)$ is even and $G(x,\cdot)$ is strictly monotone.\par
Then,  for every $\lambda>0$, there exists $\mu^*> 0$ such that, for each $\mu>\mu^*$ and for each convex set $C\subseteq L^2(\Omega)$ whose closure
in $L^2(\Omega)$ contains the set $\{G(\cdot,u(\cdot)) : u\in H^1_0(\Omega)\}$, 
there exists $v^*\in C$ such
that the problem
$$\cases {-\left ( a+b\int_{\Omega}|\nabla u(x)|^2dx\right )\Delta u =\mu f(x,u)-\lambda v^*(x)g(x,u) & in $\Omega$\cr
& \cr u=0 & on $\partial\Omega$\cr}$$
has at least three weak solutions, two of which are global minima in $H^1_0(\Omega)$ of the functional
$$u\to {{a}\over {2}}\int_{\Omega}|\nabla u(x)|^2dx+{{b}\over {4}}\left ( \int_{\Omega}|\nabla u(x)|^2dx\right ) ^2-\mu\int_{\Omega}F(x,u(x))dx
+\lambda\int_{\Omega}v^*(x)G(x,u(x))dx\ .$$}\par
\smallskip
PROOF. Set
$$D=\left \{x\in \Omega : \sup_{\xi\in {\bf R}}F(x,\xi)>0\right\}\ .$$
By assumption, meas$(D)>0$. Then, by  the Scorza-Dragoni theorem,
there exists a compact set $K\subset D$, of positive measure, such that the restriction
of $F$ to $K\times {\bf R}$ is continuous. Fix a point $\hat x\in K$ such that
the intersection of $K$ and any ball centered at $\hat x$ has a positive measure.
Choose $\hat \xi\in {\bf R}$ so that $F(\hat x,\hat\xi)>0$. By continuity, there
is $r>0$ such that
$$F(x,\hat\xi)>0$$
for all $x\in K\cap B(\hat x,r)$. Set
$$M=\sup_{(x,\xi)\in \Omega\times [-|\hat \xi|,|\hat \xi|]}|F(x,\xi)|\ .$$
Since $f\in {\cal A}$, we have $M<+\infty$. Next, choose an open set $\tilde\Omega$
such that
$$K\cap B(\hat x,r)\subset\tilde\Omega\subset\Omega$$
and
$$\hbox {\rm meas}(\tilde\Omega\setminus (K\cap B(\hat x,r)))<{{\int_{K\cap B(\hat x,r)}F(x,\hat \xi)dx}\over {M}}\ .$$
Finally, choose a function $\tilde u\in H^1_0(\Omega)$ such that
$$\tilde u(x)=\hat\xi$$
for all $x\in K\cap B(x,r)$,
$$\tilde u(x)=0$$
for all $x\in \Omega\setminus\tilde\Omega$ and
$$|\tilde u(x)|\leq |\hat\xi|$$
for all $x\in\Omega$. Thus, we have
$$\int_{\Omega}F(x,\tilde u(x))dx=\int_{K\cap B(\hat x,r)}F(x,\hat \xi)dx+\int_{\tilde\Omega\setminus (K\cap B(\hat x,r))}F(x,\tilde u(x))dx$$
$$>\int_{K\cap B(\hat x,r)}F(x,\hat \xi)dx-M\hbox {\rm meas}(\tilde\Omega\setminus (K\cap B(\hat x,r)))>0\ .$$
Now, fix any $\lambda>0$ and set
$$\mu^*={{\Phi(\tilde u)+{{\lambda}\over {2}}I_{Gg}(\tilde u)}\over {I_f(\tilde u)}}\ .$$
Fix $\mu>\mu^*$. Hence
$$\Phi(\tilde u)-\mu I_f(\tilde u)+{{\lambda}\over {2}}I_{Gg}(\tilde u)<0\ .$$
From this, we infer that the functional $\Phi-\mu I_f+{{\lambda}\over {2}}I_{Gg}$ possesses at least to global
minima since it is even. At this point, we can apply Theorem 1 to the functions $g$ and $\mu f-\lambda Gg$. Our current
conclusion follows from the one of Theorem 1
since we have $\lambda^*=0$ and hence we can take the same fixed $\lambda>0$.\hfill $\bigtriangleup$
\medskip
We conclude by proposing two very challenging problems assuming $n\geq 4$.\par
\medskip
PROBLEM 1. - Does Proposition 1 hold if $p={{n+2}\over {n-2}}$ ?\par
\medskip
PROBLEM 2. - Keeping the assumptions on $f$, does Theorem 2 hold for any $g\in \tilde{\cal A}$ ?
\bigskip
\bigskip
\bigskip
\bigskip
{\bf Acknowledgement.} The author has been supported by the Gruppo Nazionale per l'Analisi Matematica, la Probabilit\`a e 
le loro Applicazioni (GNAMPA) of the Istituto Nazionale di Alta Matematica (INdAM). \par

\vfill\eject
\centerline {\bf References}\par
\bigskip
\bigskip
\noindent
[1]\hskip 5pt C. O. ALVES, F. S. J. A. CORR\^EA and T. F. MA,
{\it Positive solutions for a quasilinear elliptic equations
of Kirchhoff type}, Comput. Math. Appl., {\bf 49} (2005), 85-93.\par
\smallskip
\noindent
[2]\hskip 5pt M. CHIPOT and B. LOVAT, {\it Some remarks on non local
elliptic and parabolic problems}, Nonlinear Anal., {\bf 30} (1997),
4619-4627.\par
\smallskip
\noindent
[3]\hskip 5pt Z. DENKOWSKI, S. MIG\'ORSKI and N. S. PAPAGEORGIOU,
{\it An Introduction to Nonlinear Analysis: Theory}, Kluwer Academic
Publishers, 2003.\par
\smallskip
\noindent
[4]\hskip 5pt X. HE and W. ZOU, {\it Infinitely many positive
solutions for Kirchhoff-type problems}, Nonlinear Anal.,
{\bf 70} (2009), 1407-1414.
\par
\smallskip
\noindent
[5]\hskip 5pt A. MAO and Z. ZHANG, {\it Sign-changing and multiple
solutions of Kirchhoff type problems without the P. S. condition},
Nonlinear Anal., {\bf 70} (2009), 1275-1287.\par
\smallskip
\noindent
[6]\hskip 5pt K. PERERA and Z. T. ZHANG, {\it Nontrivial solutions of
Kirchhoff-type problems via the Yang index}, J. Differential Equations,
{\bf 221} (2006), 246-255.\par
\smallskip
\noindent
[7]\hskip 5pt P. PUCCI and J. SERRIN, {\it A mountain pass theorem},
J. Differential Equations, {\bf 60} (1985), 142-149.\par
\smallskip
\noindent
[8]\hskip 5pt B. RICCERI, {\it A strict minimax inequality criterion and some of its consequences}, Positivity, {\bf 16} (2012), 455-470.\par
\smallskip
\noindent
[9]\hskip 5pt E. ZEIDLER, {\it Nonlinear functional analysis and its
applications}, vol. III, Springer-Verlag, 1985.\par
\smallskip
\noindent
[10]\hskip 5pt Z. T. ZHANG and K. PERERA, {\it Sign changing solutions of
Kirchhoff type problems via invariant sets of descent flow}, J. Math.
Anal. Appl., {\bf 317} (2006), 456-463.\par
\bigskip
\bigskip
\bigskip
Department of Mathematics\par
University of Catania\par
Viale A. Doria 6\par
95125 Catania\par
Italy
{\it e-mail address}: ricceri@dmi.unict.it
\bye